\documentclass[11pt]{amsart}
\usepackage{geometry}                
\usepackage{graphicx}
\usepackage{color}
\usepackage{amssymb}
\usepackage{epstopdf}
\usepackage{float}

\usepackage{amsaddr}

\newtheorem{remark}{Remark}
\newtheorem{prop}{Proposition}
\newtheorem{teo}{Theorem}

\title[A charge-preserving method for solving GRAND networks]{A charge-preserving method for solving \\graph neural diffusion networks}

\author{Lidia Aceto\textsuperscript{\S} and  Pietro Antonio Grassi\textsuperscript{\ddag}}
\address{\textsuperscript{\S,\ddag}Università degli Studi del Piemonte Orientale -- Dipartimento di Scienze e Innovazione Tecnologica, Viale T. Michel, 11 - 15121 Alessandria}
\email{lidia.aceto@uniupo.it}
\address{\textsuperscript{\ddag}INFN and Centro Regge -- v. P. Giuria, 1 - 10125 Torino}
\email{pietro.grassi@uniupo.it}

\date{\today}                                           

\begin{document}
\begin{abstract}
The aim of this paper is to give a systematic mathematical interpretation of the diffusion problem on which 
Graph Neural Networks (GNNs) models are based.
The starting point of our approach is a dissipative functional leading to dynamical equations which allows us to study the symmetries of the model. We discuss the conserved charges and provide a charge-preserving numerical method for solving the dynamical equations. In any dynamical system and also in GRAph Neural Diffusion (GRAND), knowing the charge values and their conservation along the evolution flow could provide a way to understand how GNNs and other networks work with their learning capabilities.
\end{abstract}



\maketitle

\section{Introduction} \label{sec:1}

Graph Neural Networks (GNNs) are a very powerful architecture for neural networks and deep learning algorithms. They are generally based on multiple-layer structures to encode the information, 
non-linear activation functions and suitable cost functions to compare the forwarded input with the expected results of the training data set. This essential structure has been implemented in different variants with excellent results. Nonetheless, these models have some flaws which fails to exploit the full power of the GNNs (see e.g. \cite{MB} for a comprehensive discussion). To overcome some of these problems, a very promising technique which combines ideas from graph neural networks and diffusion models 
has been explored with outstanding performances. This technique is known as GRAND, acronym for GRAph Neural Diffusion \cite{GRAND,di2022graph}.
The diffusion coefficient is the \textit{attention matrix} \cite{GAT}, which depends nonlinearly on the features of the nodes, and the evolution of the system (forwarding) is obtained by numerically solving the diffusion equation. The above-mentioned layers are substituted by the steps in the numerical approximation and the number of layers is optimized for maximum efficiency. The implementation of this new technique has been put forward in several papers \cite{MB,GRAND,GREAND,chamberlain2021beltrami}, tackling the problem by discretizing the diffusion equation in time and solving it with different numerical methods such as forward and backward Euler methods and Runge-Kutta methods, and optimizing the methods by some ad hoc prescriptions. Some of these methods are more efficient than others, and a limited number of applications have been developed for some of them. Nevertheless, the simplicity of the implementation makes the approach very relevant for the AI community. In this paper, we provide a more systematic approach both from theoretical point of view (the construction of a functional implementing the symmetries of the model -- in the spirit of \cite{di2022graph}) and from numerical implementation by 
providing a well-suited technique. Furthermore, this systematic work has been done to reveal the real knowledge of the network through the conserved charges of the dynamical system. 

In order to improve the efficiency and explainability of the neural networks, we adopt a novel approach based on physical principles by studying the symmetries of the diffusion differential equation and implementing a numerical algorithm to take these symmetries into account in the learning process. The knowledge of these symmetries is definitely of great importance as it allows a deeper understanding of the mechanism underlying the GNNs -- which contributes to the {\it explainability} of neural networks -- and the discovery that the network has learned some of the constants of motion corresponding to conserved quantities under the symmetry group.

Given the diffusion equation, there are many symmetries that preserve the form of the equation, and they are widely studied in the mathematical literature. However, among these symmetries there are also some that are unessential for our purposes. To find out the really interesting symmetries (which are respected by the data set), we use the well-known Noether's theorem. It states that for each symmetry of the theory there is a conserved quantity that contains the physical information about the system. The most powerful approach to study the symmetries in this way is to formulate the problem using an action functional, to derive the Ward identities and consequently the conserved quantities. In our case, the specific nature of the differential equation, which is a parabolic differential equation, falls outside the normal techniques of the energy (or action) functional approach (see \cite{evans2022partial} and the discussion in \cite{digiovanni2023understanding}). In this case, the differential equation cannot be derived from an energy functional, since the energy itself is not a conserved quantity being a dispersive process on the graph. A more appropriate technique will be the {\it gradient  flow} mechanism. However, the latter does not lead to a simple treatment of the symmetries. Finally, we adopt a technique proposed in \cite{ortiz,gradient-flow} where a new functional is proposed to deal with the gradient flow. The action functional corresponds to a network $\sigma$-model (where the $\sigma$-model maps are vectors 
-- named feature vectors -- from graph vertices to a $d$-dimensional vector space) coupled to a parametric time-dependent background. The Euler-Lagrangian (EL) equations derived from this functional coincide with the diffusion equation only in a suitable limit of the background, far from that limit the EL equations are of the second order and they can be easily formulated in a Hamiltonian framework by introducing the phase space. The advantage of the phase space formulation is simplified expressions for the conserved charges. In that way, we have a well-suited formulation of the  equations with the relevant symmetries manifest, starting point for the numerical analysis. 
 
 Before discussing the numerical approach, we point out that the difference between the present discussion and the case of the Schr\"odinger equation is the existence for the latter of a complex structure for the vector space, that allows a simple action functional formulation and the energy is conserved together with charges correspondent to the symmetries of the model. In our case, we do not assume any complex structure of the vector space and, for that reason we adopt the formulation of 
 \cite{gradient-flow}. 
 
Once that the equations and the conserved charges have been constructed we should wonder whether the conservation of the charges are affected by the numerical scheme used to study the evolution of the system. Unfortunately, this is not automatically implemented by any numerical scheme and therefore we should adopt the most suitable one. It turns out that the conserved charges discussed in the present work appear to be quadratic polynomials on the phase space and they are preserved by the second order numerical scheme known as {\it implicit midpoint scheme} discussed for example in \cite{hairer}. That scheme, in contrast to more used schemes such as forward or backward  Euler schemes,  preserves all quadratic invariants. Notice that our diffusion equation is non linear, nonetheless the charge expression is quadratic polynomial. However, the compatibility between the numerical scheme and the conservation of the charges is possible if the learnable parameters respect those symmetries. Those parameters are learned by the well-known backpropagation algorithm trimmed on a set of reference data. If the data, and consequently the parameter inherit those symmetries they arrange themselves to carry this information on the classification of the results, but that it is possible if the numerical scheme adopted preserves the symmetries. 

We would like to mention some related works on the subject. Recently, appear several physical inspired neural network models \cite{lutter2019deep,cranmer2020lagrangian,greydanus2019hamiltonian} based on Lagrangian and Hamiltonian dynamics. The evolution equation of the system are those obtained by the Hamilton-Jacobi equations or from Euler-Lagrangian equations. Nonetheless, none of those works tackle the problem of diffusivity in their scheme and how the conserved quantities are obtained in that case. On the other side, interesting work appeared in \cite{liu2021machine,mattheakis2020physical,bondesan2019learning,liu2022machine}, where the symmetries are learned from the data on the form of trajectories of a given dynamical system. The approach is different from ours, but it would be very interesting to see if given the trajectories derived from a dissipative equations one can learn all possible symmetries preserved by the system. That will matter of future investigations. Finally, the present work is a companion of a more applicative work in collaboration with computer science colleagues and the main goal would be to see whether the symmetries would be enough to reproduce the classification algorithm achieved in similar works. The network learns the most natural quantities, those which persist in the evolution of the system, like humans
\cite{wang2023scientific,forestano2023deep,hagemeyer2022learning}. 

The paper is organized as follows. In Section~\ref{sec:2} we establish some notations and recall preliminary notions of graph theory useful in the following. In Section~\ref{sec:3} we introduce a functional which allow us to write the dynamical system depending on a parameter $\epsilon$ which, as  $\epsilon \rightarrow  0,$ reduces to the diffusion equation already used in GRAND. 
In Section~\ref{sec:4} we discuss conserved quantities of the continuous dynamical system. Then, Section~\ref{sec:5} is devoted to the discretization of the continuous problem using classical numerical methods, including the implicit midpoint scheme, in the light of preserving the charges along the solution with the numerical methods. 
Numerical experiments 
in support of the present analysis
are reported in Section~\ref{sec:6}.
Finally, we give some concluding remarks in Section~\ref{sec:7}. In Appendix~\ref{AppA}, we discuss the model in the continuum. We compare the graph model studied here with the conventional $\sigma$-model used in theoretical physics pointing out the non-local and non-linear structures.

\section{Background} \label{sec:2}

We briefly review some basic notions of graph theory and collect some notations that will be used in the following discussions. 

An undirected (self-avoiding) graph is a pair $\mathcal{G} = (\mathcal{V}, \mathcal{E})$ where  $\mathcal{V}=\{1,2,\dots,n\}$ is a set of nodes  (or vertices) and $\mathcal{E} \subseteq \mathcal{V} \times \mathcal{V}$ is a set of edges such that
if $(i, j) \in \mathcal{E},$ then $(j, i) \in \mathcal{E},$  for all $i \neq j.$ The topology of the graph is encoded into the {\it incidence matrix} 
\begin{eqnarray}\nonumber 
   w_{ij} =\left\{ \begin{array}{ll}
         1 & {\rm if} \,\, (i,j) \in \mathcal{E} \\
         0 & {\rm otherwise}
    \end{array} \right. \,. 
\end{eqnarray}

Let $\mathbf{x}:\mathcal{V} \times [t_0,T] \longrightarrow \mathbb{R}^d$ denote a function defined on nodes and a positive time interval such that the {\it feature vector} is given by 
\[
\mathbf{x}_{i}(t) = (\begin{array}{cccc}
  \!  x_i^1(t), & x_i^2(t),& \dots, & x_i^d (t) 
\end{array} \!)^T.
\]  
We can define differential operators acting on such function. The {\it graph gradient} is the operator $\nabla$ mapping functions defined on vertices to functions defined on edges:
\begin{eqnarray}
\label{intA}
(\nabla \mathbf{x})_{ji} = \mathbf{x}_{j}(t)- \mathbf{x}_{i}(t).
\end{eqnarray}
It is immediate to check that $(\nabla \mathbf{x})_{ij} = -(\nabla \mathbf{x})_{ji}.$

For two column vectors $\mathbf{u}, \mathbf{v} \in \mathbb{R}^{d},$ such that $\mathbf{u}=(u^1,u^2,\dots,u^d)^T,$ we adopt the notation 
\[ 
\mathbf{u}^{T} \mathbf{v} = \sum_{I=1}^{d} {u}^{I} {v}^{I}
\]
for the Euclidean scalar product; hence, the Euclidean norm can be written in a coordinate-free way as
\begin{eqnarray}
    \label{intB}
\| \mathbf{u} \| = \sqrt{\mathbf{u}^{T} \mathbf{u}}.
\end{eqnarray}

We represent rotations with orthogonal matrices $\Lambda \in SO(d)$ 
or, equivalently, $\Lambda = e^\mathcal{R}$
where $\mathcal{R}= \sum_{a=1}^{d(d-1)/2} c_a \mathcal{R}_a$ is 
 a linear combination of Lie-algebra-$\mathfrak{so}(d)$ skew-symmetric matrices $\mathcal{R}_a$ \cite{gilmore}.  Obviously, the Euclidean scalar product is invariant under any rotation. 

 In addition, we need the definition of partial derivatives with respect to a vector $\mathbf{u}$ which is given by 
\begin{equation} \label{dervec2}
    \frac{\partial \mathbf{u}}{\partial \mathbf{v}} = 
    \begin{pmatrix}
         \frac{\partial u^1}{\partial v^1} &\dots & \frac{\partial u^d}{\partial v^1} \\
         \vdots & \ddots & \vdots \\
         \frac{\partial u^1}{\partial v^d} & \dots& \frac{\partial u^d}{\partial v^d} \\
    \end{pmatrix}.
   \end{equation} 
Then, for a given scalar $c$ the chain rule can be specified as follows:
\begin{equation} \label{derscalchain}
    \frac{\partial c( \mathbf{u})}{\partial \mathbf{v}} =   \frac{\partial \mathbf{u}}{\partial \mathbf{v}}  \frac{\partial c( \mathbf{u})}{\partial \mathbf{u}}. 
\end{equation}
Finally, given a generic functional in the form
\[
F[q] = \int_{t_\mathrm{i}}^{t_\mathrm{f}} \mathcal{F}(t,q(t),q'(t),q''(t),\dots, q^{(m)}(t),\dots) \, dt,
\]
with $q^{(m)}(t)$ the $m$-th derivative of $q(t),$ the functional derivative of $F[q]$ is defined as
\begin{equation} \label{derfunz}
    \frac{\delta F}{\delta q} = \sum_{m \ge 0} \frac{d^m}{dt^m} \frac{\partial F}{\partial q^{(m)}}.
\end{equation}
Suitable boundary conditions on $t_\mathrm{i}$ and $t_\mathrm{f}$ should be assigned. 

\section{Variational Principle} \label{sec:3}

In order to describe the symmetries, the conservation of charges and the energy of the system, it is useful to introduce the following functional
\begin{eqnarray} \label{fuAA}
L_{\epsilon}[\mathbf{x},V] = 
\int_{t_0}^{T}e^{-\frac{t}{\epsilon}} 
\left( \frac12\sum_{i \in \mathcal{V}} \, \Big\| \frac{ \partial \mathbf{x}_{i}}{\partial t} \Big\|^2 + 
\frac{1}{ 2\epsilon} \sum_{i,j \in \mathcal{V}}  G_{ij} \, \big\| \big(\nabla \mathbf{x}\big)_{ij}\big\|^2  + V[\mathbf{x}] \right) \, dt,
\end{eqnarray}
letting $\epsilon$ be a small enough positive number. Here $G_{ij}$ is the $(i,j)$-th entry of the symmetric matrix $G$ of order $n$ with respect to the edge $(i,j)$ and  $V$ denotes the potential that depends on $\mathbf{x}$ in a nonlinear way and does not depend on time and discrete-space derivatives.  Also $G_{ij}$ depends on $\mathbf{x}$ in a nonlinear way (see \cite{GAT} and reference therein). Before continuing, we clarify this fact by considering the form of this matrix in more detail. 

In \eqref{fuAA}  $G_{ij}$  enters in the action multiplied by $ \big\| \big(\nabla \mathbf{x}\big)_{ij}\big\|^2$ which is symmetric under the exchange of $i,j$ 
(see \eqref{intA} and \eqref{intB}). Therefore, only the symmetric part of $G_{ij}$ enters in the action and the skew-symmetric part drops out. This can be easily implemented by restring the dependence of $G_{ij}$ upon symmetric invariants. Introducing two  matrices   $W_K, W_Q$ ($K$ stands for keys, $Q$ for queries \cite{GAT, GRAND, GREAND}) which are of size $d'\times d$ with $d'\gg d$, for all $i,j$ we can construct the following expression
\begin{eqnarray} \label{cicA}
    {\mathcal M}_{ij} = ( W_K \mathbf{x}_{i})^T (W_Q \mathbf{x}_{j})\,.
\end{eqnarray}
Thus, setting
\begin{equation} \label{cij}
    C_{ij}= {\mathcal M}_{ij}+{\mathcal M}_{ji} \,,
\end{equation}
the entries of the matrix $G$ can be given as
\begin{eqnarray} \label{cicB}
    {G}_{ij} = w_{ij} \, g\left( C_{ij}\right);
\end{eqnarray}
here the factor $w_{ij}$ is the incidence matrix of the graph 
that selects which node is linked to another and the {\it activation function} $g$ could be the $\operatorname{softmax}$ function defined by
\[\operatorname{softmax}\left( C_{ij}\right) = 
    \frac{e^{C_{ij}}}{\displaystyle{\sum_{j \in \mathcal{V}} e^{C_{ij}}}},
\] see \cite{higham} for further details. When $C_{ij}$ is given as in \eqref{cij} the entries of the matrix $G$ are given in terms of ``scaled~dot''. Other possibilities  are obtained by introducing the new expressions
\begin{eqnarray*} \label{cicBA}
     \mathcal{K}_i  = 
    \| W_K \mathbf{x}_{i} \|,   \qquad 
    \mathcal{Q}_j  =  \| W_Q \mathbf{x}_{j} \|
\end{eqnarray*}
and by constructing $C_{ij}$ in terms of these building blocks
\begin{itemize}
    \item   $\frac{{\mathcal M}_{ij}}{\mathcal{K}_i \mathcal{Q}_j},$ \qquad  
    \item $\exp\left(- (\mathcal{K}_i^2 + \mathcal{Q}_j^2 -  {\mathcal M}_{ij} - {\mathcal M}_{ji})/ \varsigma^2\right),\, $ with $\varsigma$ a suitable scaling. 
\end{itemize}
The $\operatorname{softmax}$ function evaluated on these two matrices is indicated in the literature as ``cosine-similarity'' and ``exponential~kernel'', respectively \cite{chamberlain2021beltrami}. Actually, the softmax function can be substituted with any activation function such as ReLU, sigmoids or another useful distribution.

The next step is to derive the equations of motion from the functional \eqref{fuAA}. It is worth observing that  for a single graph point the second term vanishes and only the potential plays an effective role. However, although the latter is also necessary for several applications (see, for example, \cite{digiovanni2023understanding}),  to simplify the analysis and without losing generality, we will omit this term from now on and then, we use 
\begin{eqnarray} \label{fuA}
L_{\epsilon}[\mathbf{x}] := L_{\epsilon}[\mathbf{x},0] = 
\int_{t_0}^{T}e^{-\frac{t}{\epsilon}} 
\left( \frac12\sum_{i \in \mathcal{V}} \, \Big\| \frac{ \partial \mathbf{x}_{i}}{\partial t} \Big\|^2 + 
\frac{1}{ 2\epsilon} \sum_{i,j \in \mathcal{V}}  G_{ij} \, \big\| \big(\nabla \mathbf{x}\big)_{ij}\big\|^2 \right) \, dt.
\end{eqnarray}
Denoting by 
\begin{equation} \label{lagrangian}
    \mathcal{L}_\epsilon \Big(t,\mathbf{x}_{i},\frac{\partial \mathbf{x}_{i}}{\partial t}\Big)=e^{-\frac{t}{\epsilon}} 
\left( \frac12\sum_{i \in \mathcal{V}} \Big( \frac{ \partial \mathbf{x}_{i}}{\partial t} \Big)^{T} \frac{\partial \mathbf{x}_{i}}{\partial t} + 
\frac{1}{ 2\epsilon} \sum_{i,j \in \mathcal{V}}  G_{ij} \big( \big(\nabla \mathbf{x}\big)_{ij}\big)^T   \big(\nabla \mathbf{x}\big)_{ij}  \right)
\end{equation}
the Lagrangian, a local functional which can be viewed as $t$-dependent-background field theory $\sigma$-model,
we can rewrite (\ref{fuA}) as follows
\[
L_{\epsilon}[\mathbf{x}] = 
\int_{t_0}^{T} \mathcal{L}_\epsilon \Big(t,\mathbf{x}_{i},\frac{\partial \mathbf{x}_{i}}{\partial t}\Big) \, dt.
\]
Consequently, taking the functional derivatives of $L_{\epsilon}[\mathbf{x}]$ with  respect to $\mathbf{x}_{i}$ we obtain  (see (\ref{derfunz}))\footnote{In deriving the differential equations we have imposed the boundary conditions $\frac{\partial \mathbf{x}_i}{\partial t}\big|_{t=t_0} =0$ while at $t=T$ we took in account the limit for $\epsilon \rightarrow 0$.}

\begin{eqnarray} \label{pez1}
    \frac{\delta L_\epsilon}{\delta \mathbf{x}_{i}} &=&   \frac{\partial \mathcal{L}_\epsilon}{\partial \mathbf{x}_{i}} - \frac{d}{dt} \frac{\partial \mathcal{L}_\epsilon}{\partial(\partial \mathbf{x}_{i}/\partial t)} \nonumber\\
    &=&  \frac{e^{-\frac{t}{\epsilon}}}{ 2\epsilon}  \frac{\partial }{\partial \mathbf{x}_{i}}\left(  \sum_{i,j \in \mathcal{V}}  G_{ij} \big( \big(\nabla \mathbf{x}\big)_{ij}\big)^T   \big(\nabla \mathbf{x}\big)_{ij} \right) -\frac{d}{dt} \left( e^{-\frac{t}{\epsilon}}  \frac{\partial \mathbf{x}_{i}}{\partial t}\right) \nonumber\\
    &=& e^{-\frac{t}{\epsilon}} \Bigg[\frac{1}{ 2\epsilon}  \frac{\partial }{\partial \mathbf{x}_{i}}\left(  \sum_{i,j \in \mathcal{V}}  G_{ij} \big( \big(\nabla \mathbf{x}\big)_{ij}\big)^T   \big(\nabla \mathbf{x}\big)_{ij} \right) +
    \frac{1}{ \epsilon}
  \frac{\partial \mathbf{x}_{i}}{\partial t} -  \frac{\partial^{2} \mathbf{x}_{i}}{\partial t^{2}} \Bigg].
\end{eqnarray}
Using \eqref{dervec2} and \eqref{derscalchain} and the symmetry of the matrix $G$, by direct calculation we can write 
\begin{eqnarray} \label{pez2}
   && \frac{\partial }{\partial \mathbf{x}_{i}}\left(  \sum_{\ell,k \in \mathcal{V}}  G_{\ell k} \big( \big(\nabla \mathbf{x}\big)_{\ell k}\big)^T   \big(\nabla \mathbf{x}\big)_{\ell k} \right) 
   \nonumber\\
   &=&
   \sum_{\ell,k \in \mathcal{V}} \Bigg[ \frac{\partial  G_{\ell k} }{\partial \mathbf{x}_{i}}   \big( \big(\nabla \mathbf{x}\big)_{\ell k}\big)^T   \big(\nabla \mathbf{x}\big)_{\ell k}  + 
   G_{\ell k}  \frac{\partial }{\partial \mathbf{x}_{i}}\left( \big( \big(\nabla \mathbf{x}\big)_{\ell k}\big)^T   \big(\nabla \mathbf{x}\big)_{\ell k} \right)\Bigg] \nonumber\\
  \nonumber\\
  &=&
   2 \sum_{k \in \mathcal{V}} \Bigg\{\Bigg[   \frac{\partial  G_{k} }{\partial \mathbf{x}_{i}}    \mathbf{x}_{k}^T  -    \sum_{\ell \in \mathcal{V}}\frac{\partial  G_{\ell k} }{\partial \mathbf{x}_{i}}   \mathbf{x}_{\ell}^T    \Bigg]    + 2   \left( 
   G_{i}  \, \delta_{ i k}  -  
 G_{i k}   \right) I_d \Bigg\} \mathbf{x}_{k}, \nonumber
\end{eqnarray}
where we have introduced the weight of the node $k$ summing the columns of $G_{k \ell}$
\[
    G_k = \sum_{\ell \in  \mathcal{V}} G_{k \ell},
\]
$\delta_{ik}$ is the Kronecker delta and $I_d$ is the identity matrix of order $d$. 
Therefore, setting 
\begin{eqnarray}
\label{fuBAA}
A_{ij} =       \sum_{\ell \in \mathcal{V}}\frac{\partial  G_{\ell j} }{\partial \mathbf{x}_{i}}   \mathbf{x}_{\ell}^T    -\frac{\partial  G_{j} }{\partial \mathbf{x}_{i}}    \mathbf{x}_{j}^T       + 2   \left(  G_{i j} -G_{i}  \, \delta_{ i j}    \right) I_d,
\end{eqnarray}
relation \eqref{pez1} becomes
\[
 \frac{\delta L_\epsilon}{\delta \mathbf{x}_{i}} 
    = e^{-\frac{t}{\epsilon}} \Bigg[-\frac{1}{ \epsilon}  \sum_{j \in \mathcal{V}} A_{ij}  \mathbf{x}_{j} +
    \frac{1}{ \epsilon}
  \frac{\partial \mathbf{x}_{i}}{\partial t} -  \frac{\partial^{2} \mathbf{x}_{i}}{\partial t^{2}} \Bigg].
\]
Notice that for fixed $i,j,$  $A_{ij}$ is a square matrix of order $d$. Imposing 
\[
\frac{\delta L_\epsilon}{\delta \mathbf{x}_{i}} =0
\]
leads to
\begin{eqnarray}
\label{fuB}
-  \frac{\partial^{2} \mathbf{x}_{i}}{\partial t^{2}} + \frac1\epsilon \left( 
 \frac{\partial \mathbf{x}_{i}}{\partial t} -   
 \sum_{j \in \mathcal{V}} A_{ij} \mathbf{x}_{j}  \right)= 0.
\end{eqnarray}
As $\epsilon \rightarrow 0$, we recover the diffusion equation {given in \cite{GRAND} on which GRAND is based, namely}
\begin{eqnarray}
\label{fuBA}
\frac{\partial \mathbf{x}_{i}}{\partial t} -   
 \sum_{j \in \mathcal{V}} A_{ij} \mathbf{x}_{j}=0.
\end{eqnarray}
We emphasize that the main difference between \eqref{fuB} and \eqref{fuBA} is the second derivative term reminiscent of dumped harmonic oscillators. This additional term allows us to define the conserved charges as we are going to discuss. 

To convert the equation (\ref{fuB}) into a first-order system we use the rescaled momenta $\mathbf{P}_i$ defined as follows
\begin{eqnarray}
\label{fuC}
\mathbf{p}_{i}(t) = \frac{\partial \mathcal{L}_\epsilon}{ \partial (\partial \mathbf{x}_i/\partial t)} =
  e^{-\frac{t}{\epsilon}} 
 \frac{\partial \mathbf{x}_{i}}{\partial t} \equiv e^{-\frac{t}{\epsilon}} \mathbf{P}_{i}(t)\,.
\end{eqnarray}
Thus, the equation (\ref{fuB}) becomes the system 
\begin{equation} \label{fuF}
\left\{
\begin{array}{ll}
\displaystyle{\frac{\partial \mathbf{x}_{i}}{\partial t} } &=  \mathbf{P}_{i} \\
 &\\
\displaystyle{ \displaystyle{\frac{\partial \mathbf{P}_{i}}{\partial t}}}&= \displaystyle{ \frac{1}{\epsilon} \Big( \mathbf{P}_{i} -
\sum_{j \in \mathcal{V}}A_{ij} \mathbf{x}_{j} \Big)}
\end{array}\right..
\end{equation}

So far we have only focused on the dynamic of a particular node. Now we consider the equations that take into account all the nodes of the graph. At this aim we introduce the
vector $\mathbf{Y}(t) = (\mathbf{x}^T_1(t), \dots, \mathbf{x}^T_n(t), 
\mathbf{P}^T_1(t), \dots, \mathbf{P}^T_n(t))^T
\in \mathbb{R}^{2 \hat{n}},$ where $\hat{n} = n d.$ 
Under the hypothesis that $G_{ij}$ depends upon the scalar products (see \eqref{cicA} and \eqref{cij})
\begin{equation} \label{W}
   {C}_{ij} =   \mathbf{x}_{i}^T \left(W_K^T W_Q+W_Q^T W_K\right) \mathbf{x}_{j} \equiv \mathbf{x}_{i}^T \mathbb{W} \mathbf{x}_{j}  
\end{equation}
we can reorganize the term $(-\sum_{j \in \mathcal{V}}A_{ij} \mathbf{x}_{j})$ given in  \eqref{fuF} and define the matrix ${\mathcal C}(\mathbf{Y}(t))$ 
having the following block-entries 
\begin{eqnarray}
    \label{suppaA}
   \qquad  \Big({\mathcal C}(\mathbf{Y}(t))\Big)_{rs} =  \left\{ \begin{array}{ll}
       \displaystyle{ -2 I_d \sum_{\substack{ i \in \mathcal{V} \\ i\neq r}} w_{ri} g\left(\mathbf{x}_r^T \mathbb{W} \mathbf{x}_i\right)}\,, &  \mbox{if } r=s \\ 
        & \\
         w_{rs}\Big( 2 I_d 
     g\left(\mathbf{x}_r^T \mathbb{W} \mathbf{x}_s\right)
    -\mathbb{W}  
    g'\left(\mathbf{x}_r^T \mathbb{W} \mathbf{x}_s\right)
    \|\mathbf{x}_r - \mathbf{x}_s\|^2\Big)\,, & \mbox{if } r\neq s
    \end{array} \right. \,,
\end{eqnarray}
for $r,s=1,2,\dots,n$; here $g'(z)=\frac{dg}{dz}.$ Then, the system can be expressed 
in matrix form as
\begin{eqnarray}\label{fuFA}
\frac{\partial \mathbf{Y}(t)}{\partial t}
&=&
\mathcal{B}(\mathbf{Y}(t)) \, \mathbf{Y}(t),  \qquad \mathcal{B}(\mathbf{Y}(t)) =  \frac{1}{\epsilon}\left(\begin{array}{cc}
       O  & \epsilon I_{\hat{n}} \\
        {\mathcal C}(\mathbf{Y}(t)) & 
        I_{
       \hat{n}}
    \end{array}\right), \quad t\in (t_0,T],
\end{eqnarray}
where $I_{\hat{n}}$ is the identity matrix  of order $\hat{n}.$ Notice that $\mathcal{C}(\mathbf{Y}(t))$ is a symmetric matrix since its entries are function of $\mathbb{W}$  which is itself a symmetric matrix (see \eqref{W}). In the subsequent analysis this property will be crucial.

\section{Conserved Quantities}\label{sec:4}

For any dynamical system, it is useful to study possible constants of motion, associated to conserved quantities, which may characterize the evolution in a simpler way. The most natural example is the energy of the system which is the conserved quantity associated to the Hamiltonian. Furthermore, there might be other conserved quantities associate to the symmetries of the model, such as rotational symmetry, scale symmetry or translation symmetry. If the Hamiltonian exhibits one of those symmetries, there are such conserved quantities. 

In this case the Hamiltonian is (see \eqref{lagrangian} and \eqref{fuC}) 
\begin{eqnarray}
\label{fuD}
H \big( \mathbf{x}_{i}, \mathbf{P}_{i}\big) &=& \sum_{i \in \mathcal{V}} e^{-\frac{t}{\epsilon}} \mathbf{P}^{T}_{i}  \frac{\partial \mathbf{x}_{i}}{\partial t} - \mathcal{L}_\epsilon \nonumber\\
&=& \frac{e^{-\frac{t}{\epsilon}}}{2} \sum_{i \in \mathcal{V}} \Big( \frac{ \partial \mathbf{x}_{i}}{\partial t} \Big)^{T}   \frac{\partial \mathbf{x}_{i}}{\partial t} 
-\frac{e^{-\frac{t}{\epsilon}}}{ 2\epsilon} \sum_{i,j \in \mathcal{V}}  G_{ij} \big(\nabla \mathbf{x}\big)_{ij}^T  \big(\nabla \mathbf{x}\big)_{ij} \nonumber\\
&=&
 e^{-\frac{t}{\epsilon}}
\left(    \frac{1}{2}  \sum_{i \in \mathcal{V}} \mathbf{P}^{T}_{i} \mathbf{P}_{i} -
 \frac{1}{ 2\epsilon} \sum_{i,j \in \mathcal{V}}  G_{ij}  \big(\nabla \mathbf{x}\big)_{ij}^T  \big(\nabla \mathbf{x}\big)_{ij} \right)
\end{eqnarray}
and, as a consequence of  the explicit dependence upon $t$ in \eqref{fuD}, we have 
\begin{eqnarray}
\label{fuH}
\frac{d H}{d t} &=& 
 \frac{e^{-\frac{t}{\epsilon}}}{\epsilon} 
 \sum_{i \in \mathcal{V}} \mathbf{P}^{T}_{i}\mathbf{P}_{i} - \frac1\epsilon H 
 \nonumber \\
&=& 
\frac{e^{-\frac{t}{\epsilon}} }{2\epsilon}  
\left(      \sum_{i \in \mathcal{V}} \mathbf{P}^{T}_{i} \mathbf{P}_{i} +
 \frac{1}{\epsilon} \sum_{i,j \in \mathcal{V}}  G_{ij} \big(\nabla \mathbf{x}\big)_{ij}^T  \big(\nabla \mathbf{x}\big)_{ij} \right),
\end{eqnarray}
where in the first equality we used \eqref{fuF}. 
The Hamiltonian  is conserved in the limit $\epsilon~\rightarrow~0$ (which corresponds also to $T\rightarrow \infty$) 
namely when the dispersive functional $L_{\epsilon}[\mathbf{x}]$ is independent of $t$ and the energy is a well-defined concept. Therefore, the energy is not a valuable quantity in the present case. 

On the other hand, we study other constants of motion that are conserved along the diffusion process and can be used to analyze the network final result.  For that, we notice that while  the expressions
$\mathbf{P}^{T}_{i} \mathbf{P}_{i}$ and 
$ \big(\nabla \mathbf{x}\big)_{ij}^T  \big(\nabla \mathbf{x}\big)_{ij}  $  in \eqref{fuD}
are invariant under any rotation, the 
invariance of $G_{ij}$ is restricted by the presence of 
$\mathbb{W}.$
Denoting by  
\[\Lambda = e^{\mathcal R} = I_d + {\mathcal R} + {\mathcal O}({\mathcal R}^2), \qquad \Lambda \in SO(d),\] 
a generic rotation -- where ${\mathcal R}$ is a skew-symmetric matrix -- those 
which are preserved commute with $\mathbb{W}$, namely
 \begin{eqnarray}
     \label{WR}
     {[\mathbb{W}, \mathcal{R}] = O}\,. 
 \end{eqnarray}
In the case where  $G_{ij}$ depends upon ${\mathcal M}_{ij},$  $ \mathcal{K}_i,$ and  $\mathcal{Q}_j,$ then the preserved rotations are those which commute with $\mathbb{W}, W^T_K W_K,$ and $W^T_Q W_Q$. Obviously, that subgroup is smaller than the one given by \eqref{WR}. The matrices $W_K$ and $W_Q$ are the learnable parameters of the theory and are obtained using the backpropagation technique. Therefore, after a few iterations of the machine learning process, they stabilized to some specific form and the symmetry should be compatible with those.

Since the Hamiltonian given in (\ref{fuD}) is invariant under the following rigid symmetry 
\begin{eqnarray}
\label{syA}
\mathbf{x}_{i}  \,\, \leftarrow \,\, \Lambda  \mathbf{x}_{i}\,,  \qquad \mathbf{P}_{i}  \,\,  \leftarrow  \,\,  \Lambda  \mathbf{P}_{i},
\end{eqnarray}
we can define 
the corresponding  charge as 
\begin{eqnarray}
\label{syB}
Q_\Lambda (t)= e^{-\frac{t}{\epsilon}} \sum_{i \in \mathcal{V}} 
\left(\mathbf{P}_{i}(t)\right)^{T} {\mathcal R}\mathbf{x}_{i}(t)\,. 
\end{eqnarray}
Using $\mathbf{Y}(t)$, we can rewrite this equation as 
\begin{eqnarray}
    \label{syBA}
  Q_\Lambda(t) = -  \frac{e^{-\frac{t}{\epsilon}}}{2}
  \left(\mathbf{Y}(t)\right)^{T}
   (\mathcal{J}  \otimes \mathcal{R}) \mathbf{Y}(t)
\end{eqnarray}
where $\mathcal{J}$ is the skew-symmetric matrix
\[
\mathcal{J} = \left( \begin{array}{rr}
     O  & I_n \\
     -I_n  & O
  \end{array}
    \right).
\]
By a direct calculation we have that
\begin{eqnarray}
\label{syC}
\frac{d Q_\Lambda}{dt}  &=&
- \frac{1}{\epsilon} Q_\Lambda - 
 \frac{e^{-\frac{t}{\epsilon}}}{2} 
 \frac{d }{dt} \big(\mathbf{Y}^T  (\mathcal{J}  \otimes \mathcal{R}) \mathbf{Y} \big)
 \nonumber \\
 &=&
 \frac{e^{-\frac{t}{\epsilon}}}{2\epsilon}  \bigg(\mathbf{Y}^T  (\mathcal{J}  \otimes \mathcal{R}) \mathbf{Y} \bigg) - 
 \frac{e^{-\frac{t}{\epsilon}}}{2} 
 \frac{d }{dt} \big(\mathbf{Y}^T  (\mathcal{J}  \otimes \mathcal{R}) \mathbf{Y} \big)
 \nonumber \\
 &=&
  \frac{e^{-\frac{t}{\epsilon}}}{2}  \bigg(\frac{1}{\epsilon}\mathbf{Y}^T  (\mathcal{J}  \otimes \mathcal{R}) \mathbf{Y}  - 
 \frac{d }{dt} \big(\mathbf{Y}^T  (\mathcal{J}  \otimes \mathcal{R}) \mathbf{Y} \big) \bigg)
 \nonumber \\
 &=&
  \frac{e^{-\frac{t}{\epsilon}}}{2}  \bigg(\frac{1}{\epsilon}\mathbf{Y}^T  (\mathcal{J}  \otimes \mathcal{R}) \mathbf{Y}  - \mathbf{Y}^T 
  \Big(\mathcal{B}(\mathbf{Y})^T \mathcal{J} \otimes \mathcal{R}  +   \mathcal{J} \otimes \mathcal{R} \,
\mathcal{B}(\mathbf{Y})\Big)
  \mathbf{Y}.
  \bigg) 
\end{eqnarray}
Therefore 
\begin{eqnarray} \label{conserva}
\frac{d Q_\Lambda}{dt} =0 \quad &\Longleftrightarrow& \quad 
\mathbf{Y}^T    
\Bigg(
\frac{1}{\epsilon}  (\mathcal{J}  \otimes \mathcal{R})-
\mathcal{B}(\mathbf{Y})^T \mathcal{J} \otimes \mathcal{R}  +   \mathcal{J} \otimes \mathcal{R} 
\mathcal{B}(\mathbf{Y})\Bigg)\mathbf{Y} =0 \,.
\end{eqnarray}
When $[\mathbb{W}, \mathcal{R}]=O,$ which expresses the compatibility of the learnable parameters with the simmetries, from a straightforward computation we can check that
\begin{equation} \label{relmain}
\Big(\mathcal{B}(\mathbf{Y})^T \mathcal{J} \otimes \mathcal{R}  +   \mathcal{J} \otimes \mathcal{R} \,
\mathcal{B}(\mathbf{Y})\Big) =  \frac1\epsilon  \mathcal{J} \otimes \mathcal{R}.  
\end{equation}
and this implies that  $Q_\Lambda(t)$ actually is a constant.
We summarize this fact in the following result.
\begin{teo}
  Let $\mathbb{W}$ be the matrix defined in \eqref{W} and $\mathcal{R}$ any skew-symmetric matrix of order $d.$ With reference to \eqref{fuFA}, when
\[
[\mathbb{W}, \mathcal{R}]=O,
\]
the charge $Q_\Lambda$ given in \eqref{syBA} is preserved, that is 
\[
Q_\Lambda(t) = Q_\Lambda(t_0), 
\qquad \forall \, t \in (t_0,T].
\] 
\end{teo}

Since we need to solve numerically the dynamical system \eqref{fuFA}, in the next section we are interested in selecting an appropriate numerical scheme that is 
able to inherit the conservation of the underlying continuous invariants. In particular, in this analysis we are interested to select a method that preserves the charge.

\section{Computing discrete invariants} \label{sec:5}

Instead of discretizing the dynamical system \eqref{fuFA}, from a numerical point of view it is more convenient to consider the system in the canonical coordinates   $(\mathbf{x}_{i}, \mathbf{p}_{i}),$ 
(see \eqref{fuC} and \eqref{fuF}),
which in matrix form  can be represented as 
\begin{eqnarray}\label{xpmat}
\frac{\partial \mathbf{y}(t)}{\partial t}
&=&
\mathcal{E}(\mathbf{y}(t)) \, \mathbf{y}(t),  \qquad \mathcal{E}(\mathbf{y}(t)) =  \left(\begin{array}{cc}
       O  & e^{\frac{t}{\epsilon}} I_{\hat{n}} \\
       -\frac{e^{-\frac{t}{\epsilon}}}{\epsilon}{\mathcal C}(\mathbf{y}(t)) & 
        O
    \end{array}\right), \quad t\in (t_0,T],
\end{eqnarray}
with $\mathbf{y}(t) = (\mathbf{x}^T_1(t), \dots, \mathbf{x}^T_n(t), 
\mathbf{p}^T_1(t), \dots, \mathbf{p}^T_n(t))^T.$ In this framework the charge is given by (see \eqref{syB}) 
\begin{eqnarray} \label{syBa}
Q_\Lambda (t) =   \sum_{i \in \mathcal{V}} 
\left(\mathbf{p}_{i}(t)\right)^{T} {\mathcal R}\mathbf{x}_{i}(t)  = 
-\frac{1}{2}  \left(\mathbf{y}(t)\right)^{T} (\mathcal{J}  \otimes \mathcal{R}) \mathbf{y}(t).
\end{eqnarray}
Consequently,
\begin{equation*}  
\frac{d Q_\Lambda}{dt} =0 \quad \Longleftrightarrow \quad   
\mathbf{y}^T    
\Bigg( \mathcal{E}(\mathbf{y})^T \mathcal{J} \otimes \mathcal{R}  +   \mathcal{J} \otimes \mathcal{R} \,
\mathcal{E}(\mathbf{y})\Bigg)\mathbf{y} =0 \,.
\end{equation*} 
Under the hypothesis that $[\mathbb{W}, \mathcal{R}]=O,$ from a direct computation we can check that
\begin{equation} \label{eqprinc}
    \Big(\mathcal{E}(\mathbf{y})^T \mathcal{J} \otimes \mathcal{R}  +   \mathcal{J} \otimes \mathcal{R} \,
\mathcal{E}(\mathbf{y})\Big) =  O
\end{equation}
and then the charge is trivially preserved.

Denoting by $\{\mathbf{y}^{(k)}\}_{k=0}^N$ the numerical solution provided by the  method used to discretize \eqref{xpmat}, the conservation of charge \eqref{syBa} in the discrete frame means that at two consecutive grid points the following relation should be examined
\begin{eqnarray} \label{pofD}
 (\mathbf{y}^{(k+1)})^T  (\mathcal{J}  \otimes \mathcal{R}) \mathbf{y}^{(k+1)}
  = (\mathbf{y}^{(k)})^T  (\mathcal{J}  \otimes \mathcal{R}) \mathbf{y}^{(k)}.
\end{eqnarray}

We first consider the most commonly used schemes, namely the forward Euler  and backward Euler  methods. Since both are unable to preserve the charges -- as will be proved below -- we introduce another numerical method capable of doing so. Following Hairer {\it et al.} \cite[Sec. IV.2]{hairer}, we consider the simplest case of Gauss collocation methods. In all cases, starting from an initial guess $\mathbf{y}^{(0)}= \mathbf{y}(t_0),$
we solve the equation in \eqref{xpmat} on  the uniform grid $\{t_k = t_0+k h; k=0, 1, \dots,N; t_N=T\}.$ 
\subsubsection{The forward Euler method}

The numerical solution obtained by applying the forward Euler method can be computed by
\begin{eqnarray}
    \label{pofC}
    \mathbf{y}^{(k+1)}  &=& \mathbf{y}^{(k)} + h \mathcal{E}(\mathbf{y}^{(k)})  \mathbf{y}^{(k)} 
\end{eqnarray}
where $\mathbf{y}^{(k)} \approx \mathbf{y}(t_k).$ Using this recurrence relation we can prove the following result.

\begin{prop} \label{prop1}
The forward Euler method does not preserve the charge \eqref{syBa}.
\end{prop}

\proof By substituting in \eqref{pofD} the expression of $\mathbf{y}^{(k+1)}$ provided by \eqref{pofC} 
we can write   
\begin{eqnarray}
    \label{pofE}
&& (\mathbf{y}^{(k+1)})^T  (\mathcal{J}  \otimes \mathcal{R}) \mathbf{y}^{(k+1)}
  =
  (\mathbf{y}^{(k)})^T (\mathcal{J}  \otimes \mathcal{R}) \mathbf{y}^{(k)} +\nonumber \\
  &&\hspace{1cm} + \, h (\mathbf{y}^{(k)})^T \Bigg[ \left( \mathcal{E}(\mathbf{y}^{(k)})\right)^T  (\mathcal{J}  \otimes \mathcal{R})  + (\mathcal{J}  \otimes \mathcal{R}) 
   \mathcal{E}(\mathbf{y}^{(k)}) \Bigg]\mathbf{y}^{(k)}+
  \nonumber\\
   &&\hspace{1cm} + \, h^2 (\mathbf{y}^{(k)})^T \left(\mathcal{E}(\mathbf{y}^{(k)})\right)^T  (\mathcal{J}  \otimes \mathcal{R})  \mathcal{E}(\mathbf{y}^{(k)}) \mathbf{y}^{(k)}.
\end{eqnarray}
Then, taking into account \eqref{eqprinc}
and observing that 
\begin{equation} \label{matimport}
    \left(\mathcal{E}(\mathbf{y}^{(k)})\right)^T  (\mathcal{J}  \otimes \mathcal{R})  \mathcal{E}(\mathbf{y}^{(k)}) = \frac{1}{\epsilon}
\left(\begin{array}{cc}
       O  & {\mathcal C}(\mathbf{y}^{(k)}) \\
      {\mathcal C}(\mathbf{y}^{(k)}) & 
        O
    \end{array}\right) (\mathcal{J}  \otimes \mathcal{R})
\end{equation}
we have
\[
  (\mathbf{y}^{(k+1)})^T  (\mathcal{J}  \otimes \mathcal{R}) \mathbf{y}^{(k+1)}
  =
  (\mathbf{y}^{(k)})^T (\mathcal{J}  \otimes \mathcal{R}) \mathbf{y}^{(k)} + \mathcal{O}\left(\frac{h^2}{\epsilon} \right).
\]
This last relation leads us to the conclusion that the forward Euler  method does not preserve the charge (see \eqref{pofD}).
\endproof
As just proved, the conservation of the charge is spoiled at order $h^2$. 
Furthermore, we observe that for $\epsilon \rightarrow \infty$ the charge is conserved. This is expected since in that limit we recover a simple model which conserve any charge. On the other side, in the limit 
$\epsilon \rightarrow 0$, that is when the diffusion equation is recovered, the second term could become arbitrarily large. Thus, we have to choose a different discretization scheme. 

\subsubsection{The backward Euler method}
When we discretize the system \eqref{xpmat} by the backward Euler method we get
\begin{eqnarray*}
    \mathbf{y}^{(k+1)}  &=& \mathbf{y}^{(k)} + h \mathcal{E}(\mathbf{y}^{(k+1)})  \mathbf{y}^{(k+1)}
\end{eqnarray*}
and then
\begin{eqnarray}
    \label{BEsol}
   \left(I_{2 \hat{n}} - h \mathcal{E}(\mathbf{y}^{(k+1)}) \right) \mathbf{y}^{(k+1)}  &=& \mathbf{y}^{(k)}. 
\end{eqnarray}

Using this recurrence relation we can prove the following result.
\begin{prop}
The backward Euler method does not preserve the charge \eqref{syBa}.
\end{prop}

\proof By substituting in \eqref{pofD} the expression of $\mathbf{y}^{(k)}$ provided by \eqref{BEsol}, we can write   
\begin{eqnarray*}
 (\mathbf{y}^{(k)})^T  (\mathcal{J}  \otimes \mathcal{R}) \mathbf{y}^{(k)}
  &=&
  (\mathbf{y}^{(k+1)})^T  \left(I_{2 \hat{n}} - h \mathcal{E}(\mathbf{y}^{(k+1)}) \right)^T (\mathcal{J}  \otimes \mathcal{R})  \left(I_{2 \hat{n}} - h \mathcal{E}(\mathbf{y}^{(k+1)}) \right)\mathbf{y}^{(k+1)} \nonumber\\
 &=&   (\mathbf{y}^{(k+1)})^T   (\mathcal{J}  \otimes \mathcal{R}) \mathbf{y}^{(k+1)} \nonumber \\
   &-&  h (\mathbf{y}^{(k+1)})^T \Bigg[ \left( \mathcal{E}(\mathbf{y}^{(k+1)})\right)^T  (\mathcal{J}  \otimes \mathcal{R})  + (\mathcal{J}  \otimes \mathcal{R}) 
   \mathcal{E}(\mathbf{y}^{(k+1)}) \Bigg]\mathbf{y}^{(k+1)}\nonumber\\
     &+&  h^2 (\mathbf{y}^{(k+1)})^T \Bigg[ \left( \mathcal{E}(\mathbf{y}^{(k+1)})\right)^T (\mathcal{J}  \otimes \mathcal{R}) \mathcal{E}(\mathbf{y}^{(k+1)}) \Bigg]\mathbf{y}^{(k+1)}.
\end{eqnarray*}
Using \eqref{eqprinc} and considering the expression of the matrix in \eqref{matimport} we obtain 
 \begin{eqnarray}
    \label{BEcharge}
(\mathbf{y}^{(k)})^T  (\mathcal{J}  \otimes \mathcal{R}) \mathbf{y}^{(k)}
  &=&
   (\mathbf{y}^{(k+1)})^T  (\mathcal{J}  \otimes \mathcal{R}) \mathbf{y}^{(k+1)}+
   \mathcal{O}\left(\frac{h^2}{\epsilon} \right).
\end{eqnarray}
and therefore  the charge is not conserved.
\endproof

\subsubsection{The implicit midpoint method}
The scheme we propose to use in this setting works as follows 
\begin{eqnarray} \label{pofH}
    \mathbf{y}^{(k+1)}  &=& \mathbf{y}^{(k)} + h \mathcal{E}(\mathbf{y}^{(\xi)})  \frac{\mathbf{y}^{(k+1)}+\mathbf{y}^{(k)}}{2}, 
\end{eqnarray}
where $\mathbf{y}^{(\xi)}$ approximates the continuous solution at a point $t_\xi \in [t_n, t_{n+1}]$.  Recombining the above expression gives
\[
\left( I_{2 \hat{n}} -\frac{h}{2} \mathcal{E}(\mathbf{y}^{(\xi)})\right)  \mathbf{y}^{(k+1)}  = \left( I_{2 \hat{n}} +\frac{h}{2} \mathcal{E}(\mathbf{y}^{(\xi)})\right)\mathbf{y}^{(k)}.
\]
Since (see \eqref{xpmat})
\[
 I_{2 \hat{n}} - \frac{h}{2} \mathcal{E}(\mathbf{y}^{(\xi)}) = 
 \left(\begin{array}{cc}
       I_{\hat{n}}  & - \frac{h \, e^{\frac{t}{\epsilon}}}{2} I_{\hat{n}} \\
       \frac{h \, e^{-\frac{t}{\epsilon}}}{2\epsilon}{\mathcal C}(\mathbf{y}^{(\xi)}) & 
        I_{\hat{n}}
    \end{array}\right),
\]
the blocks on the main diagonal are both invertible and then this matrix can be inverted blockwise as follows
\[
\begin{split}
 \left( I_{2 \hat{n}} -\frac{h}{2} \mathcal{E}(\mathbf{y}^{(\xi)}) \right)^{-1}= &\,
 \Bigg( I_2 \otimes \Big[ I_{\hat{n}} + \frac{h^2}{4\epsilon}   {\mathcal C}(\mathbf{y}^{(\xi)}) \Big]^{-1} \Bigg)
\left( I_{2 \hat{n}} +\frac{h}{2} \mathcal{E}(\mathbf{y}^{(\xi)}) \right).
\end{split}
\]
From the above considerations and by a direct calculation we  have
\begin{eqnarray} \label{mid}
    \mathbf{y}^{(k+1)}  &= &\,  \Bigg( I_2 \otimes \Big[ I_{\hat{n}} + \frac{h^2}{4\epsilon}   {\mathcal C}(\mathbf{y}^{(\xi)}) \Big]^{-1} \Bigg)\left( I_{2 \hat{n}} +\frac{h}{2} \mathcal{E}(\mathbf{y}^{(\xi)})\right)^2\mathbf{y}^{(k)} \nonumber \\
    &= &\, \left( I_{2 \hat{n}} +\frac{h}{2} \mathcal{E}(\mathbf{y}^{(\xi)})\right)^2  \Bigg( I_2 \otimes \Big[ I_{\hat{n}} + \frac{h^2}{4\epsilon}   {\mathcal C}(\mathbf{y}^{(\xi)}) \Big]^{-1} \Bigg)\mathbf{y}^{(k)}  \nonumber\\
   & \equiv & \, \mathcal{P}(\mathbf{y}^{(\xi)}) \mathbf{y}^{(k)}.
\end{eqnarray}
Now, denoting for the reader's convenience 
\[
\begin{split}
\mathcal{G}(\mathbf{y}^{(\xi)}) :=  
& \,   \left(I_{\hat{n}}  + \frac{h^2}{4\epsilon}  {\mathcal C}(\mathbf{y}^{(\xi)}) \right)^2  \left( I_{n} \otimes \mathcal{R} \right)
 \end{split}
\]
and recalling that  
\[
-{\mathcal C}(\mathbf{y}^{(\xi)}) \,(I_n\otimes \mathcal{R}) +(I_n\otimes \mathcal{R}) \,{\mathcal C}(\mathbf{y}^{(\xi)}) =O
\]
since $(I_n\otimes \mathcal{R})$ is a block diagonal matrix whose entries commute with $\mathbb{W},$ we get 
\[
\begin{split}
&\, \left( I_{2 \hat{n}} +\frac{h}{2} \left(\mathcal{E}(\mathbf{y}^{(\xi)}) \right)^T\right)^2 (\mathcal{J}  \otimes \mathcal{R}) 
     \left( I_{2 \hat{n}} +\frac{h}{2} \mathcal{E}(\mathbf{y}^{(\xi)})\right)^2 =  \left(\begin{array}{cc}
       O & \mathcal{G}(\mathbf{y}^{(\xi)}) \\
      \left(\mathcal{G}(\mathbf{y}^{(\xi)}) \right)^T & 
        O
        \end{array}\right) \\
      =  &\, \left(\begin{array}{cc}
       O & \left(I_{\hat{n}}  + \frac{h^2}{4\epsilon}  {\mathcal C}(\mathbf{y}^{(\xi)}) \right)^2  \left( I_{n} \otimes \mathcal{R} \right) \\
     - \left(I_{\hat{n}}  + \frac{h^2}{4\epsilon}  {\mathcal C}(\mathbf{y}^{(\xi)}) \right)^2  \left( I_{n} \otimes \mathcal{R} \right)& 
        O
        \end{array}\right) \,;
        \end{split}
\]
then 
\[
\left( \mathcal{P}(\mathbf{y}^{(\xi)})\right)^T (\mathcal{J}  \otimes \mathcal{R})   \mathcal{P}(\mathbf{y}^{(\xi)}) = \mathcal{J}  \otimes \mathcal{R}.
\]
This equality together with \eqref{mid} allows us to write
\begin{eqnarray*}
 (\mathbf{y}^{(k+1)})^T  (\mathcal{J}  \otimes \mathcal{R}) \mathbf{y}^{(k+1)}
  =
   (\mathbf{y}^{(k)})^T    (\mathcal{J}  \otimes \mathcal{R}) \mathbf{y}^{(k)}.
\end{eqnarray*}
\begin{remark} \label{remark1}
It is important to emphasize that any point  $t_\xi$ is suitable for charge conservation. In particular, setting $t_\xi = t_k + \frac{h}{2} \equiv t_{k+\frac{1}{2}},$ the method \eqref{pofH} is an $A$-stable implicit method of the second order known in the literature as {\rm implicit midpoint method} (the simple instance of Gauss collocation methods, with one collocation point per step).
Actually this method can be seen as a sequence of backward Euler and then forward Euler methods. In fact, applying both methods step by step on a mesh with stepsize $h/2$ we have
\begin{eqnarray*}
   \mathbf{y}^{(k+\frac{1}{2})} &=& \mathbf{y}^{(k)} + \frac{h}{2} \mathcal{E}(\mathbf{y}^{(k+\frac{1}{2})})  \mathbf{y}^{(k+\frac{1}{2})} \\
   \mathbf{y}^{(k+1)} &=& \mathbf{y}^{(k+\frac{1}{2})} + \frac{h}{2} \mathcal{E}(\mathbf{y}^{(k+\frac{1}{2})})  \mathbf{y}^{(k+\frac{1}{2})}
\end{eqnarray*}
from which we deduce that 
\[
\mathbf{y}^{(k+\frac{1}{2})} - \mathbf{y}^{(k)}  =  \mathbf{y}^{(k+1)} -\mathbf{y}^{(k+\frac{1}{2})},
\]
or, equivalently,
\[
\mathbf{y}^{(k+\frac{1}{2})}   =  \frac{\mathbf{y}^{(k+1)} + \mathbf{y}^{(k)}}{2}.
\]
\end{remark}
We are now in the position to enunciate the following proposition. 
\begin{prop}
    The implicit midpoint method  applyied for solving the system \eqref{xpmat} reads as
\begin{eqnarray} \label{midpoint}
    \mathbf{y}^{(k+1)}  &=& \mathbf{y}^{(k)} + h \mathcal{E}\left(\frac{\mathbf{y}^{(k+1)}+\mathbf{y}^{(k)}}{2}\right)  \frac{\mathbf{y}^{(k+1)}+\mathbf{y}^{(k)}}{2}.
\end{eqnarray}
    It preserves the charge \eqref{syBa} which is a quadratic integral invariant for the system.
\end{prop}
 
\section{Numerical experiments} \label{sec:6}

In this section we report the results of some numerical experiments we have performed to show in a simple case the choices we have made in this work.

We consider a graph with three nodes (i.e. $n=3$) with vector features of dimension $d=4,$ whose dynamic is modeled by the continuous system of the form~\eqref{xpmat} defined on the interval $(0, 1].$ The incidence matrix is defined as
\[
w=
\left(\begin{array}{rrr}
 0 & 1& 1 \\
 1 & 0& 1\\
 1 & 1 &0
 \end{array} \right)
 \]
so that each node is connected to all others. As the activation function $g,$ we use the softmax function. 
Here we focus on two tests corresponding to two different ways of setting the matrix $\mathbb{W}.$ In the first test we set
$\mathbb{W} =  10^{-3} I_4,$ while in the second one we consider $\mathbb{W} = \operatorname{diag}(10^{-3},10^{-3},1,1).$ We will later refer to these two matrices as $\mathbb{W}_1$ and $\mathbb{W}_2,$ respectively. To solve the corresponding two problems, we use the forward Euler method (FE) \eqref{pofC} and the (modified) implicit midpoint method (IM), that is the method given in \eqref{pofH} with the matrix $\mathcal{E}(\mathbf{y}^{(\xi)})$ evaluated at $t_\xi = t_k.$ Although the latter method converges with order one, the numerical approximation of the solution to the next grid point can be achieved by solving only a linear system per step rather than a nonlinear problem as would occur using the implicit midpoint method. This obviously leads to a significant saving on computational cost. Furthermore, as already highlighted in the previous section, this modification does not alter the system's charge conservation. 

For the initial value $\mathbf{y}^{(0)} \in \mathbb{R}^{24}$ with 
\[
\begin{split}
\mathbf{x}_1^{(0)} & \, =(0,1,1,1)^T, \\ \mathbf{x}_j^{(0)} & \, =(1,1,1,1)^T,\quad  j=2,3, \\\mathbf{p}_\ell^{(0)} & \, =(0,0,0,0)^T, \quad \ell=1,2,3,   \\
\end{split}
\]
we calculate the numerical solution on a uniform grid with step size $h=1/50.$ In this case, the charge associated to the rotation $\mathcal{R}$ which commutes with $\mathbb{W}$ is $Q_\Lambda=0.$ 

One of the rotation matrices that we considered in our tests is the matrix
\[
\mathcal{R}_{2,1} = 
\left(\begin{array}{rrrr}
        0 & -1& 0& 0 \\
       1 & 0& 0& 0 \\
        0 & 0 &0& 0\\
        0 & 0 &0& 0\\
     \end{array} \right).
 \]
The charge here called $Q_1$ is associated with this matrix. In a similar way we can define the other rotations we used. In Table~\ref{tab:1} we name the charge associated with each of these rotations.
\begin{table}[H] 
    \centering
    \renewcommand{\arraystretch}{1.5} 
    \begin{tabular}{|c|cccccc|}
     \hline 
      $\mathcal{R}$   & $\mathcal{R}_{2,1}$  & $\mathcal{R}_{3,1}$  & $\mathcal{R}_{4,1}$  &  $\mathcal{R}_{3,2}$ &  $\mathcal{R}_{4,2}$ & $\mathcal{R}_{4,3}$ \\
      \hline
       $Q_\Lambda$  & $Q_1$  & $Q_2$ &  $Q_3$ & $Q_4$ &  $Q_5$& $Q_6$ \\
        \hline 
    \end{tabular}
    \vspace{2mm}
    \caption{Correspondence rotation-charge.}
    \label{tab:1}
\end{table}

\begin{figure}[!h]
\centering
\includegraphics[width=14.5cm,height=20cm,keepaspectratio]{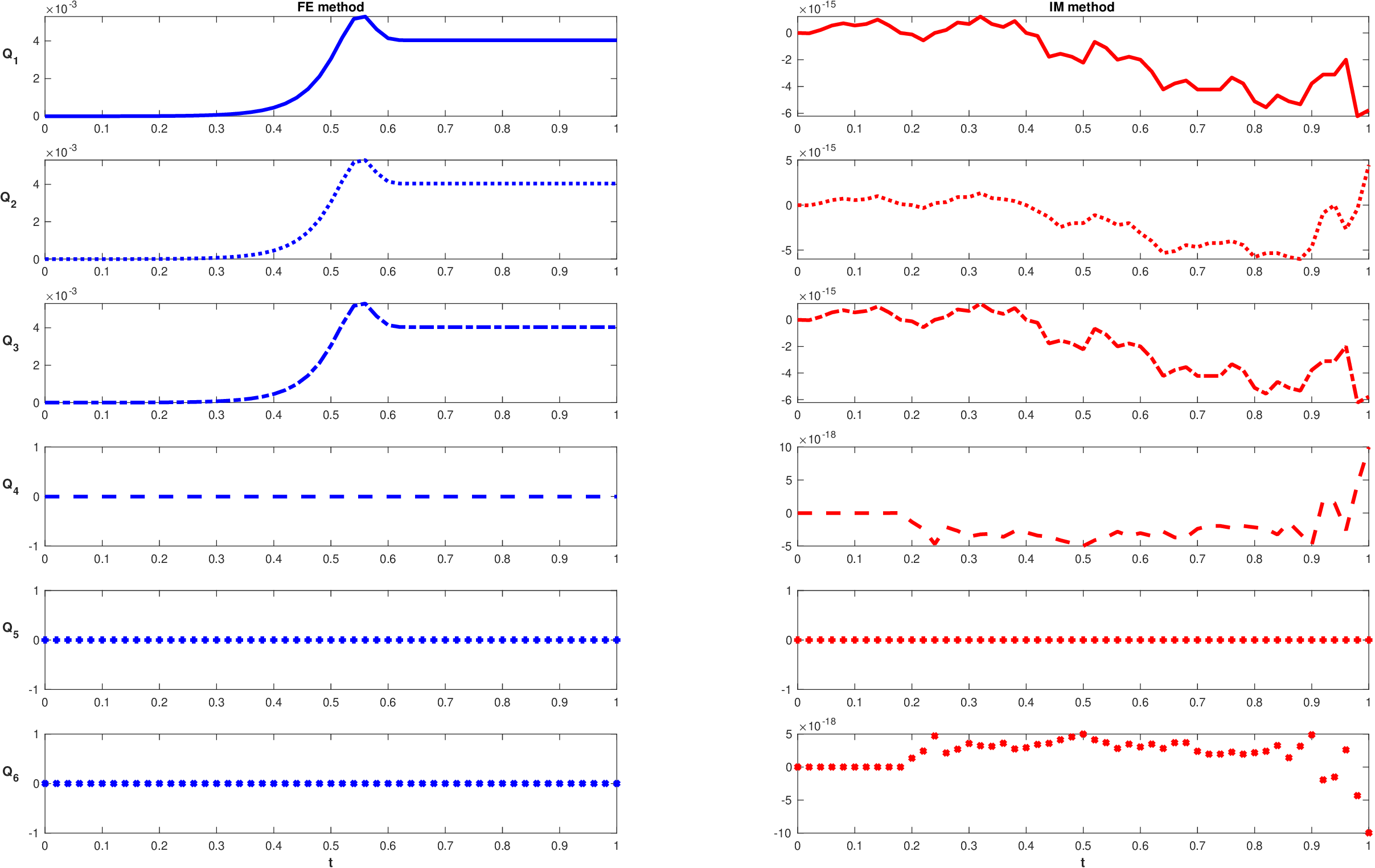}
\caption{Charge of the system ~\eqref{xpmat} with $\epsilon = 0.1$ and $\mathbb{W}_1$ from $Q_1$ up to $Q_6$ on the left computed by the forward Euler (FE) method and on the right by the (modified) implicit midpoint  (IM) method.}\label{fig:1}
\end{figure}

Since $\mathbb{W}_1$ commutes with all considered rotations, we expect according to the theoretical results of the previous section that the corresponding charges are preserved by the implicit midpoint method at each grid point; however, we generally do not expect this in the forward Euler method. This fact is confirmed by the images on the right and left in Figure~\ref{fig:1}, where we plot the charges from $Q_1$ to $Q_6$ against time $t_k=kh, k=0,1,\dots,50$ on both sides.

\begin{figure}[H]
\centering
\includegraphics[width=14.5cm,height=20cm,keepaspectratio]{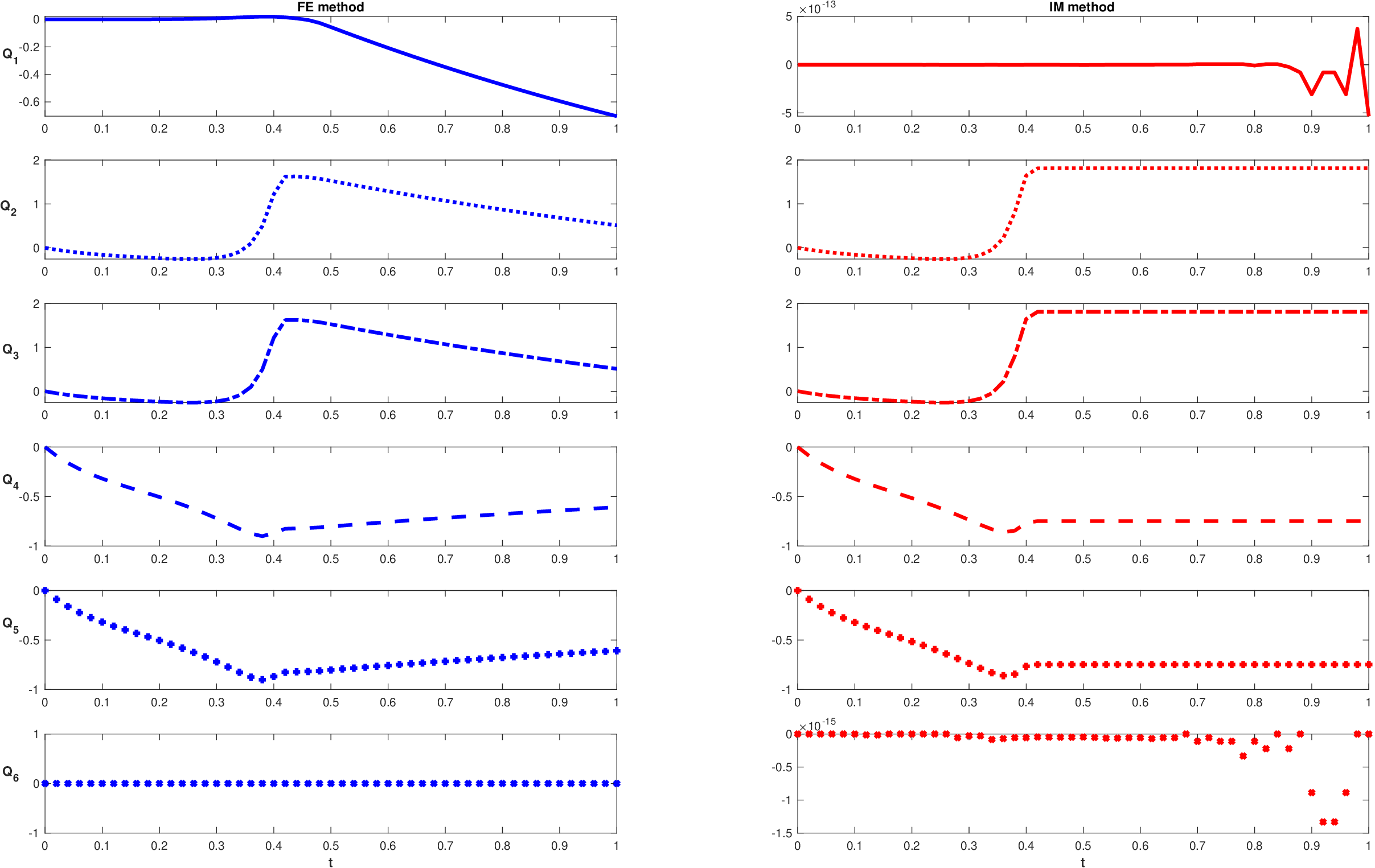}
\caption{Charge of the system ~\eqref{xpmat} with $\epsilon = 0.1$ and $\mathbb{W}_2$ from $Q_1$ up to $Q_6$ on the left computed by the forward Euler (FE) method and on the right by the (modified) implicit midpoint  (IM) method.}\label{fig:2}
\end{figure}

On the other hand, since $\mathbb{W}_2$ commutes only with $\mathcal{R}_{2,1}$ and $\mathcal{R}_{4,3}$, we expect the charges $Q_1$ and $Q_6$ to be preserved by the implicit midpoint method.
This result is confirmed by the first and last images on the right-hand side in Figure~\ref{fig:2}.
As can be seen on the left side of the same figure, this does not happen when using the forward Euler method.

In both figures, we have given the results for $\epsilon = 0.1.$ However, it should be emphasised that similar results can be obtained for other values of this parameter.

\section{Conclusions} \label{sec:7}

We systematically analysed the mathematical framework of some GNNs based on differential equations. In particular, we studied the example of the diffusion equation on a graph. We translated the study of the equation into a functional language that exploits all the relevant symmetries of the model, and we provided the best-suited numerical scheme for the system of equations.
The present work aimed to provide a unifying framework for GNNs that could have several potential applications. Complementary to more efficiency and implementability, we believe that this framework could bring a small step forward in understanding how GNNs capture the essential features from the data by learning the {\it constants of motions}, namely the conserved quantities of dynamic evolution, as the human brain does. 

\section*{Acknowledgements}
We thank M. Dossena, C. Irwin, S. Montani and L. Portinale for valuable discussions. P.A.G. would like to thank M. Caselle for encouraging him to explore these new paths.  


\appendix


\section{A model on a smooth surface} \label{AppA}

In this appendix, we provide a continuum version of the model considered in the text. 
Instead of the graph, we consider a 
$\mathbb{R}^n \times [t_0,T]$  
worldvolume  and the fields $\vec{\phi}(x,t)$ 
where $x,t$ are the coordinates (we replaced the nodes $i \in \mathcal{V}$ 
with local coordinates $x$). The features vectors ${\bf x}_{i}(t)$ on the nodes are 
replaced by continuum fields $\vec{\phi}(x,t)$ which are vectors in the feature space 
$\mathbb{R}^{d}$ as well. 

The action is 
\begin{eqnarray}
\label{supA}\nonumber
I_{\epsilon}[\vec{\phi},V] = 
\int_{t_0}^{T}  dt   
e^{-\frac{t}{\epsilon}} 
\left( \int d^nx \frac12 \frac{\partial \vec{\phi}^{T}}{\partial t} 
\frac{\partial \vec{\phi}}{\partial t} + \frac{1}{2\epsilon} \int d^nx d^ny
 \nabla^x_\mu \vec{\phi}^{T} K[x,y] 
 \nabla^{y,\mu} \vec{\phi} +
 V[\vec{\phi}]\right) 
\end{eqnarray}
where $K[x,y]$ is a non-trivial kernel and $\nabla^x_{\mu} \vec{\phi}, \nabla^{y}_{\mu} \vec{\phi}$ are the derivatives with respect to the coordinates $x,y$. The second term is 
written in terms of a second integration over $y$'s. This is due to the non-linear non-local nature of the 
model as $G_{ij}$ in \eqref{fuA}. 
We added to the Lagrangian a potential $V[\vec{\phi}]$ in order 
to describe further interactions. 

Compared with 
the graph theory, the kernel is 
$$K[x,y] = K[\vec{\phi}(x,t),\vec{\phi}(y,t)]$$
namely a function of the field $\vec{\phi}(x,t)$. Introducing the function
$$
C(x,y) = \vec{\phi}^{T}(x,t) \mathbb{W} \vec{\phi}(y,t)\,, 
$$
we can construct the bilocal expression $K[\vec{\phi}(x,t),\vec{\phi}(y,t)]$ 
as a function of $C(x,y).$
For example, a possible 
non-linear non-local expression used in 
GNNs is
\begin{eqnarray}
\label{supAB}\nonumber
K[\vec{\phi}(x,t),\vec{\phi}(y,t)] = 
\frac{ e^{C(x,y)}}
{\int d^{n}y 
 e^{C(x,y)}}\,, ~~~~~~~~~
\end{eqnarray}
with some additional conditions to make the expression well-defined. 

Computing the functional derivative with respect to $\vec{\phi}(z,t)$ we get the equation of motion
\begin{eqnarray}
    \label{supB}
    - \frac{\partial^2 \phi^I(z,t)}{\partial t^2} &+& \frac1\epsilon 
    \left[\frac{\partial \phi^I(z,t)}{\partial t} - \nabla_{{\phi}}^I V + \int d^ny 
    \Big(\left. \nabla^x_\mu K[x,z]\right|_{x=y}   +\right. \nonumber \\
&+& \left. 
    \nabla^x_\mu K[z,x]\right|_{x=y}\Big) \nabla^{y,\mu} {\phi}^I +
    \int d^nx d^ny 
 \nabla^x_\mu \vec{\phi}^{T} \frac{\partial K[x,y]}{\partial \phi^I(z,t)} 
 \nabla^{y,\mu} \vec{\phi} 
    \Big] =0 \,.
\end{eqnarray}
In the limit $\epsilon \rightarrow 0$, we retrieve the Gradient Flow equation 
\begin{eqnarray}
    \label{supC}
     \frac{\partial \phi^I}{\partial t} - \nabla_{{\phi}}^I V  &+&  \int d^ny 
    \Big(\left. \nabla^x_\mu K[x,z]\right|_{x=y}  + \left. 
    \nabla^x_\mu K[z,x]\right|_{x=y}\Big) \nabla^{y,\mu} {\phi}^I +
     \nonumber \\
&+& \int d^nx d^ny 
 \nabla^x_\mu \vec{\phi}^{T} \, \frac{\partial K[x,y]} {\partial \phi^I(z,t)} 
 \, \nabla^{y,\mu} \vec{\phi} 
    =0.
\end{eqnarray}
If $K[x,y] = \delta(x-y)$, \eqref{supC} reduces to heat equation.


\end{document}